\documentclass[a4paper,10pt,reqno]{amsart}
\usepackage[a4paper,tmargin=40mm,bmargin=35mm,hmargin=30mm]{geometry}
\usepackage[T1]{fontenc}
\usepackage{lmodern}
\usepackage{amsmath}
\usepackage{amssymb}
\usepackage{amsthm}
\usepackage{stmaryrd}
\usepackage{tikz}
\usetikzlibrary{cd}
\usepackage{booktabs}
\usepackage{multirow}
\usepackage{multicol}
\usepackage{here}
\usepackage{aliascnt}
\usepackage{hyperref}
\usepackage[noabbrev]{cleveref}

\newcommand{\NewTheorem}[2]{
	\newaliascnt{#1}{TheoremEnvironment}
	\newtheorem{#1}[#1]{#1}
	\aliascntresetthe{#1}
	\crefname{#1}{#1}{#2}
	\Crefname{#1}{#1}{#2}
}

\theoremstyle{definition}

\NewTheorem{Definition}{Definitions}
\NewTheorem{Remark}{Remarks}
\NewTheorem{Axiom}{Axioms}
\NewTheorem{Example}{Examples}
\NewTheorem{Observation}{Observations}
\NewTheorem{Convention}{Conventions}
\NewTheorem{Notation}{Notations}
\NewTheorem{Setting}{Settings}
\NewTheorem{Question}{Questions}
\NewTheorem{Answer}{Answers}
\NewTheorem{Conjecture}{Conjectures}
\NewTheorem{Problem}{Problems}
\NewTheorem{Solution}{Solutions}
\NewTheorem{Goal}{Goals}
\NewTheorem{Comment}{Comments}
\NewTheorem{Aim}{Aims}
\NewTheorem{Caution}{Cautions}
\NewTheorem{Exercise}{Exercises}

\theoremstyle{plain}
\NewTheorem{Proposition}{Propositions}
\NewTheorem{Lemma}{Lemmas}
\NewTheorem{Theorem}{Theorems}
\NewTheorem{Corollary}{Corollaries}

\crefname{enumi}{}{}
\Crefname{enumi}{}{}
\creflabelformat{enumi}{(#2#1#3)}
\crefname{enumii}{}{}
\Crefname{enumii}{}{}
\creflabelformat{enumii}{(#2#1#3)}
\crefname{enumiii}{}{}
\Crefname{enumiii}{}{}
\creflabelformat{enumiii}{(#2#1#3)}

\makeatletter
\renewcommand{\p@enumii}{}
\renewcommand{\p@enumiii}{}
\makeatother

\numberwithin{equation}{section}
\crefname{equation}{}{}
\Crefname{equation}{}{}
\creflabelformat{equation}{(#2#1#3)}

\newcommand{\SwapSymbols}[1]{
	\expandafter\let\expandafter\temporarysymbol\csname #1\endcsname
	\expandafter\let\csname #1\expandafter\endcsname\csname var#1\endcsname
	\expandafter\let\csname var#1\endcsname\temporarysymbol
}

\SwapSymbols{epsilon}
\SwapSymbols{phi}
\SwapSymbols{Gamma}
\SwapSymbols{Delta}
\SwapSymbols{Theta}
\SwapSymbols{Lambda}
\SwapSymbols{Xi}
\SwapSymbols{Pi}
\SwapSymbols{Sigma}
\SwapSymbols{Upsilon}
\SwapSymbols{Phi}
\SwapSymbols{Psi}
\SwapSymbols{Omega}

\newcommand{\cA}{\mathcal{A}}
\newcommand{\cB}{\mathcal{B}}
\newcommand{\cC}{\mathcal{C}}

\newcommand{\cG}{\mathcal{G}}

\newcommand{\cX}{\mathcal{X}}

\newcommand{\kp}{\mathfrak{p}}
\newcommand{\kq}{\mathfrak{q}}

\newcommand{\set}[2][]{\mathopen{#1\{}#2\mathclose{#1\}}}
\newcommand{\setwithspace}[2][]{\mathopen{#1\{}\,#2\,\mathclose{#1\}}}

\newcommand{\setwithcondition}[3][]{\mathopen{#1\{}\,#2\mathrel{#1|}#3\,\mathclose{#1\}}}

\DeclareMathOperator{\Hom}{Hom}
\DeclareMathOperator{\End}{End}

\DeclareMathOperator{\Mod}{Mod}
\let\mod\relax
\DeclareMathOperator{\mod}{mod}

\DeclareMathOperator{\noeth}{noeth}
\DeclareMathOperator{\QCoh}{QCoh}

\DeclareMathOperator{\Ann}{Ann}

\DeclareMathOperator{\Ker}{Ker}

\let\Im\relax
\DeclareMathOperator{\Im}{Im}

\DeclareMathOperator{\Spec}{Spec}

\DeclareMathOperator{\Min}{Min}
\DeclareMathOperator{\Ass}{Ass}
\DeclareMathOperator{\Supp}{Supp}

\DeclareMathOperator{\height}{ht}

\DeclareMathOperator{\ASpec}{ASpec}

\DeclareMathOperator{\AMin}{AMin}
\DeclareMathOperator{\AAss}{AAss}
\DeclareMathOperator{\ASupp}{ASupp}

\title{Finiteness of the number of minimal atoms in Grothendieck categories}

\author{Ryo Kanda}
\address{Department of Mathematics, Graduate School of Science, Osaka University, Toyonaka, Osaka, 560-0043, Japan}
\email{ryo.kanda.math@gmail.com}
\begin{document}

%%%%%%%%%%%%%%%%%%%%%%%%%%%%%%%%%%%%%%%%%%%%%%%%%%%%%%%%%%%%%%%%%%%%%%%%%%%%%%%%
\begin{abstract}
	For a Grothendieck category having a noetherian generator, we prove that there are only finitely many minimal atoms. This is a noncommutative analogue of the fact that every noetherian scheme has only finitely many irreducible components. It is also shown that each minimal atom is represented by a compressible object.
\end{abstract}
%%%%%%%%%%%%%%%%%%%%%%%%%%%%%%%%%%%%%%%%%%%%%%%%%%%%%%%%%%%%%%%%%%%%%%%%%%%%%%%%

\maketitle
\tableofcontents

%%%%%%%%%%%%%%%%%%%%%%%%%%%%%%%%%%%%%%%%%%%%%%%%%%%%%%%%%%%%%%%%%%%%%%%%%%%%%%%%
\section{Introduction}
\label{sec.Intro}
%%%%%%%%%%%%%%%%%%%%%%%%%%%%%%%%%%%%%%%%%%%%%%%%%%%%%%%%%%%%%%%%%%%%%%%%%%%%%%%%

The Gabriel spectrum \cite{MR0232821} of a Grothendieck category is the set of isomorphism classes of indecomposable injective objects. It is a fundamental tool to study Grothendieck categories and has been investigated by many authors (\cite{MR1434441, MR1426488, MR1347919, MR1899866} for example). One of its significance is that it plays a role analogous to the (not necessarily closed) points of a scheme. Indeed, if a scheme $X$ is locally noetherian, then the Gabriel spectrum of the category $\QCoh X$ of quasi-coherent sheaves on $X$ is canonically identified with the set of all points of $X$.

In this viewpoint, our main result is an abelian-categorical (or noncommutative) analogue of the fact that every noetherian scheme has only finitely many irreducible components. Before giving the precise statement, we recall the notion of atom spectrum.

The \emph{atom spectrum} of an abelian category is a variant of the Gabriel spectrum, which consists of \emph{atoms} introduced by Storrer \cite{MR0360717}. The atom spectrum is defined to be the set of equivalence classes of monoform objects so that it works for abelian categories that do not have enough injective objects. For a locally noetherian Grothendieck category $\cA$, such as the category $\Mod\Lambda$ of right modules over a right noetherian ring $\Lambda$ or $\QCoh X$ for a noetherian scheme $X$, there is a canonical bijection between the atom spectrum and the Gabriel spectrum (\cite[Theorem~5.9]{MR2964615}). It is known that the localizing subcategories of $\cA$ are classified using the Gabriel spectrum by Herzog \cite{MR1434441} and Krause \cite{MR1426488}, and an analogous correspondence for a noetherian abelian category was explicitly given in \cite{MR2964615} (see \cref{ClassifOfLocSubcatAndSerreSubcat}).

The atom spectrum always has the structure of a partially ordered set, and for a commutative ring $R$, the atom spectrum of $\Mod R$ is isomorphic to the prime spectrum $\Spec R$ regarded as a partially ordered set by the inclusion. Similarly, for a locally noetherian scheme $X$, the partial order of the atom spectrum of $\QCoh X$ coincides with the specialization order, that is, $x\leq y$ if and only if $y$ belongs to the closure of $x$. Therefore the minimal atoms with respect to the partial order are analogous to the irreducible components of the scheme. If a Grothendieck category has a noetherian generator, then there exist enough minimal atoms in the sense that each atom is greater than or equal to one of the minimal atoms (\cref{ExistenceOfMinimalAtom}). Our main result is stated as follows:

\begin{Theorem}[\cref{AMinOfGrothCatWithNoethGenIsFinite} and \cref{AMinOfNoethAbelCatWithGenIsFinite}]\label{IntroAMinIsFinite}
	Let $\cA$ be either a Grothendieck category having a noetherian generator or a noetherian abelian category having a generator. Then there are only finitely many minimal atoms in $\cA$.
\end{Theorem}

The main result is a generalization of the well-known fact that every commutative noetherian ring has only finitely many minimal prime ideals, but the proof is quite different because the minimal atoms are not necessarily \emph{associated atoms} of the fixed noetherian generator (\cref{ProofForCommRing}). We shall develop a method using the notion of a compressible object, which has been studied, for example, in \cite{MR0164991, MR0393118, MR0360715, MR0352177, MR588447}. A key step in the proof is \cref{MinimalAtomIsRepresentedByCompressibleObj} asserting that every minimal atom is represented by a compressible object.

We regard \cref{IntroAMinIsFinite} as an analogous result to the finiteness of the number of irreducible components of a noetherian scheme $X$ since this property is easily reduced to that of each open affine subscheme. However, it should be noted that we cannot apply \cref{IntroAMinIsFinite} directly to $\QCoh X$ since it might not have a noetherian generator.

%%%%%%%%%%%%%%%%%%%%%%%%%%%%%%%%%%%%%%%%%%%%%%%%%%%%%%%%%%%%%%%%%%%%%%%%%%%%%%%%
\subsection*{Acknowledgments}
\label{subsec.Acknowledgments}
%%%%%%%%%%%%%%%%%%%%%%%%%%%%%%%%%%%%%%%%%%%%%%%%%%%%%%%%%%%%%%%%%%%%%%%%%%%%%%%%

The author would like to express his deep gratitude to Osamu Iyama for his encouragement and guidance in Nagoya University. The author thanks the anonymous referee for their valuable comments.

The author was a JSPS Overseas Research Fellow. This work was supported by JSPS KAKENHI Grant Numbers JP17K14164, JP16H06337, and JP13J00249.

%%%%%%%%%%%%%%%%%%%%%%%%%%%%%%%%%%%%%%%%%%%%%%%%%%%%%%%%%%%%%%%%%%%%%%%%%%%%%%%%
\section{Preliminaries}
\label{sec.Prelim}
%%%%%%%%%%%%%%%%%%%%%%%%%%%%%%%%%%%%%%%%%%%%%%%%%%%%%%%%%%%%%%%%%%%%%%%%%%%%%%%%

This section is devoted to preliminary results that we use in the proof of the main result. We refer the reader to \cite{MR0340375} for the general theory of abelian categories and Grothendieck categories. In order to avoid confusion, we first recall the definition of a generator.

\begin{Definition}\label{GenSetAndGen}
	Let $\cC$ be a category.
	\begin{enumerate}
		\item\label{GenSet} A \emph{generating set} is a family $\set{U_{i}}_{i\in I}$ of objects in $\cC$ such that for any objects $M,N\in\cC$ and for any $f,g\in\Hom_{\cC}(M,N)$, there exist an index $i\in I$ and $h\in\Hom_{\cC}(U_{i},M)$ satisfying $fh\neq gh$.
		\item\label{Gen} A \emph{generator} is an object $U$ in $\cC$ such that the singleton $\set{U}$ is a generating set.
	\end{enumerate}
\end{Definition}

By a slight abuse of notation, we say that a generating set is \emph{small} if the index set is small.\footnote{We fix a Grothendieck universe and refer its elements as \emph{small sets}. Every category consists of \emph{sets} of objects and morphisms, and the set $\Hom_{\cA}(M,N)$ for objects $M,N\in\cA$ is supposed to belong to the fixed universe.}

For an abelian category $\cA$ with small direct sums (meaning small coproducts), a family $\set{U_{i}}_{i\in I}$ of objects indexed by a small set $I$ is a generating set if and only if every object in $\cA$ is a quotient object of a direct sum of copies of $U_{i}$'s since we can form the canonical morphism $\bigoplus_{i\in I}(\bigoplus_{h\in\Hom_{\cA}(U_{i},M)}U_{i})\to M$ for each $M\in\cA$. Note that the existence of a generating set is one of the axioms of a Grothendieck category.

An abelian category $\cA$ is said to be \emph{noetherian} if the set of objects is in bijection with a small set and all the objects are noetherian. In this case, a family $\set{U_{i}}_{i\in I}$ of objects indexed by a small set $I$ is a generating set if and only if every object in $\cA$ is a quotient object of a finite direct sum of copies of $U_{i}$'s. Obviously $\cA$ itself, regarded as a family of objects, is a generating set of $\cA$.

From now on, let $\cA$ be an abelian category unless otherwise specified. We recall the definition of the atom spectrum of $\cA$ and some basic results. See \cite[section~3]{MR3452186} for more details.

\begin{Definition}\label{MonoformObjAndAtomEquivAndASpec}\leavevmode
	\begin{enumerate}
		\item\label{MonoformObj} A nonzero object $H$ in $\cA$ is called \emph{monoform} if for every nonzero subobject $L$ of $H$, no nonzero subobject of $H$ is isomorphic to a subobject of $H/L$.
		\item\label{AtomEquiv} For monoform objects $H_{1}$ and $H_{2}$ in $\cA$, we say that $H_{1}$ is \emph{atom-equivalent} to $H_{2}$ if there exists a nonzero subobject of $H_{1}$ that is isomorphic to a subobject of $H_{2}$. This is an equivalence relation on the class of monoform objects in $\cA$ (\cite[Proposition~2.8]{MR2964615}).
		\item\label{ASpec} The \emph{atom spectrum} $\ASpec\cA$ of $\cA$ is the quotient set of the set of monoform objects in $\cA$ by atom equivalence. Each element of $\ASpec\cA$ is called an \emph{atom} in $\cA$. For each monoform object $H$ in $\cA$, the equivalence class of $H$ is denoted by $\overline{H}$.
	\end{enumerate}
\end{Definition}

If $\cA$ has a small generating set, then $\ASpec\cA$ is in bijection with a small set. This can be proved similarly to \cite[Proposition~2.7 (2)]{MR3351569} by extending \cite[Proposition~IV.6.6]{MR0389953} to a small generating set.

If $\cA$ is a locally noetherian Grothendieck category, then there is a bijection from $\ASpec\cA$ to the Gabriel spectrum given by $\overline{H}\mapsto E(H)$, where $E(H)$ is the injective envelope of $H$ and it does not depend on the choice of the representative $H$ (see \cite[Theorem~5.9]{MR2964615}). The proof of \cite[Theorem~2.6]{MR3272068} shows that for every morphism $f\colon E(H)\to E(H)$ that is not an isomorphism, we have $f(H)=0$.

Let $\Lambda$ be a right noetherian ring and let $\alpha\in\ASpec(\Mod\Lambda)$. Then $\alpha$ is represented by a monoform (and hence uniform) right $\Lambda$-module $H$. Since $H$ has a nonzero cyclic submodule, we can assume that $H$ itself is cyclic. The argument in the previous paragraph shows that $H$, regarded as a subset of the left $\End_{\Lambda}(E(H))$-module $E(H)$, is annihilated by the unique maximal two-sided ideal of $\End_{\Lambda}(E(H))$, which consists of all endomorphisms of $E(H)$ that are not automorphisms.

For a commutative ring $R$, we have a bijection between the prime spectrum $\Spec R$ and the atom spectrum $\ASpec(\Mod R)$ given by $\kp\mapsto\overline{R/\kp}$ (\cite[p.\ 631]{MR0360717}; see also \cite[Proposition~7.2 (1)]{MR2964615}).

Moreover, we can define the support of each object as follows:

\begin{Definition}\label{ASupp}
	Let $M$ be an object in $\cA$. The \emph{atom support} of $M$ is
	\begin{equation*}
		\ASupp M:=\setwithcondition{\alpha\in\ASpec\cA}{\text{$\alpha=\overline{H}$ for some monoform subquotient $H$ of $M$}}.
	\end{equation*}
\end{Definition}

If $0\to L\to M\to N\to 0$ is a short exact sequence in $\cA$, then $\ASupp M=\ASupp L\cup\ASupp N$ holds (\cite[Proposition~3.3]{MR2964615}). For a family of objects $\set{M_{i}}_{i\in I}$ in $\cA$ indexed by a small set $I$, we have $\ASupp(\bigoplus_{i\in I}M_{i})=\bigcup_{i\in I}\ASupp M_{i}$ (\cite[Proposition~2.12 (1)]{MR3351569}).

\begin{Proposition}\label{ASuppOfGenerator}
	If $\set{U_{i}}_{i\in I}$ is a generating set in $\cA$, then $\ASpec\cA=\bigcup_{i\in I}\ASupp U_{i}$.
\end{Proposition}

\begin{proof}
	Let $H$ be a monoform object in $\cA$. Then there exist $i\in I$ and a nonzero morphism $f\colon U_{i}\to H$. Since $\Im f$ is a nonzero subobject of $H$, we have $\overline{H}=\overline{\Im f}\in\ASupp(\Im f)\subset\ASupp U_{i}$. Hence the assertion holds.
\end{proof}

The atom spectrum $\ASpec\cA$ is regarded as a topological space with respect to the \emph{localizing topology}, that is, $\cB:=\setwithcondition{\ASupp M}{M\in\cA}$ is an open basis for $\ASpec\cA$ (\cite[Proposition~3.2]{MR3351569}). Indeed, if $M_{1},M_{2}\in\cA$ and $\alpha\in\ASupp M_{1}\cap\ASupp M_{2}$, then there exists a monoform subquotient $H_{i}$ of $M_{i}$ that represents $\alpha$ for each $i=1,2$. Since $H_{1}$ and $H_{2}$ are atom-equivalent to each other, there exists a nonzero subobject $H$ of $H_{1}$ that is isomorphic to a subobject of $H_{2}$. It then follows that
\begin{equation*}
	\alpha\in\ASupp H\subset\ASupp M_{1}\cap\ASupp M_{2}.
\end{equation*}
This shows that $\ASpec\cA$ has a unique topology such that $\cB$ is an open basis for $\ASpec\cA$.

For a commutative ring $R$, a subset of $\ASpec(\Mod R)$ is open if and only if the corresponding subset $\Phi$ of $\Spec R$ is closed under specialization, that is, the conditions $\kp\in\Phi$ and $\kp\subset\kq$ imply $\kq\in\Phi$.

The localizing topology can be used to classify localizing subcategories and Serre subcategories. Recall that a \emph{Serre subcategory} of $\cA$ is a full subcategory closed under subobjects, quotient objects, and extensions. A Serre subcategory $\cX$ of a Grothendieck category $\cA$ is called \emph{localizing} if it is also closed under small direct sums.

For a full subcategory $\cX$ of $\cA$, define $\ASupp\cX:=\bigcup_{M\in\cX}\ASupp M$, which is an open subset of $\ASpec\cA$. For a subset $\Phi$ of $\ASpec\cA$, define the Serre subcategory $\ASupp^{-1}\Phi$ of $\cA$ by
\begin{equation*}
	\ASupp^{-1}\Phi:=\setwithcondition{M\in\cA}{\ASupp M\subset\Phi}.
\end{equation*}

\begin{Theorem}[{\cite[Theorem~3.8]{MR1434441}, \cite[Corollary~4.3]{MR1426488}, and \cite[Theorem~5.5]{MR2964615}}]\label{ClassifOfLocSubcatAndSerreSubcat}\leavevmode
	\begin{enumerate}
		\item\label{ClassifOfLocSubcat} Let $\cA$ be a locally noetherian Grothendieck category. Then we have a bijection
		\begin{equation*}
			\setwithspace{\text{localizing subcategories of $\cA$}}\to\setwithspace{\text{open subsets of $\ASpec\cA$}}
		\end{equation*}
		given by $\cX\mapsto\ASupp\cX$, and its inverse is given by $\Phi\mapsto\ASupp^{-1}\Phi$.
		\item\label{ClassifOfSerreSubcat} Let $\cA$ be a noetherian abelian category. Then we have a bijection
		\begin{equation*}
			\setwithspace{\text{Serre subcategories of $\cA$}}\to\setwithspace{\text{open subsets of $\ASpec\cA$}}
		\end{equation*}
		given by $\cX\mapsto\ASupp\cX$, and its inverse is given by $\Phi\mapsto\ASupp^{-1}\Phi$.
	\end{enumerate}
\end{Theorem}

\begin{proof}
	See also \cite[Theorem~6.8 and Proposition~6.5]{MR3452186} for a detailed proof of \cref{ClassifOfLocSubcat}.
\end{proof}

The \emph{specialization order} $\leq$ on $\ASpec\cA$ is defined by
\begin{equation*}
	\alpha\leq\beta\iff\alpha\in\overline{\set{\beta}}
\end{equation*}
where $\overline{\set{\beta}}$ is the closure of the singleton $\set{\beta}$ with respect to the localizing topology. The relation $\leq$ is a partial order due to the fact that $\ASpec\cA$ is a Kolmogorov space (\cite[Proposition~3.5 and section~4]{MR3351569}). The specialization order is a generalization of the inclusion relation of the prime ideals of a ring (\cite[Proposition~4.3]{MR3351569}) and of the specialization order of the points of a locally noetherian scheme (\cite[Corollary~7.7 (4)]{MR3452186}).

The set of minimal atoms in $\cA$ with respect to the specialization order, which is the main object to study in this paper, is denoted by $\AMin\cA$.

For a Serre subcategory $\cX$ of an abelian category $\cA$, the quotient category $\cA/\cX$ is again an abelian category and there is a canonical (covariant) functor $\cA\to\cA/\cX$, which is dense and exact. If $\cA$ is moreover a Grothendieck category and $\cX$ is a localizing subcategory, then $\cA/\cX$ is a Grothendieck category and the canonical functor has a fully faithful right adjoint $\cA/\cX\to\cA$ (see \cite[section~5]{MR3351569}).

Some properties of abelian categories are inherited by their quotient categories:

\begin{Proposition}\label{ImageOfGeneratorAndNoethObj}
	Let $\cX$ be a Serre subcategory of $\cA$, and let $F\colon\cA\to\cA/\cX$ be the canonical functor.
	\begin{enumerate}
		\item\label{ImageOfNoethObj} For every noetherian object $M$ in $\cA$, the object $F(M)$ in $\cA/\cX$ is noetherian.
		\item\label{ImageOfGenerator} Assume that $\cA$ is a Grothendieck category and $\cX$ is a localizing subcategory. If $\set{U_{i}}_{i\in I}$ is a generating set in $\cA$, then $\set{F(U_{i})}_{i\in I}$ is a generating set in $\cA/\cX$.
	\end{enumerate}
\end{Proposition}

\begin{proof}
	\cref{ImageOfNoethObj} Let $L'_{0}\subset L'_{1}\subset\cdots$ be an ascending chain of subobjects of $F(M)$. For each $i\geq 0$, by \cite[Corollary~4.3.10]{MR0340375}, there exists a subobject $N_{i}$ of $M$ such that $F(N_{i})=L'_{i}$ as a subobject of $F(M)$. Let $L_{i}:=\sum_{j=0}^{i}N_{j}$. Since $M$ is noetherian, the ascending chain $L_{0}\subset L_{1}\subset\cdots$ eventually stabilizes. As a subobject of $M$,
	\begin{equation*}
		F(L_{i})=F\Bigg(\sum_{j=0}^{i}N_{j}\Bigg)=\sum_{j=0}^{i}F(N_{j})=\sum_{j=0}^{i}L'_{j}=L'_{i}.
	\end{equation*}
	Hence the ascending chain $L'_{0}\subset L'_{1}\subset\cdots$ also eventually stabilizes. Therefore $F(M)$ is noetherian.
	
	\cref{ImageOfGenerator} This follows from the fact that $F$ is dense and exact and it preserves small direct sums.
\end{proof}

For each $\alpha\in\ASpec\cA$, we define the \emph{localization} $\cA_{\alpha}$ at $\alpha$ to be the quotient category of $\cA$ by the Serre subcategory
\begin{equation*}
	\cX(\alpha):=\ASupp^{-1}(\ASpec\cA\setminus\overline{\set{\alpha}}).
\end{equation*}
The canonical functor $\cA\to\cA_{\alpha}$ is denoted by $(-)_{\alpha}$. If $\cA$ is a Grothendieck category, then $\cX(\alpha)$ is a localizing subcategory of $\cA$. This is a generalization of the localization of a commutative ring at a prime ideal (\cite[Proposition~6.9]{MR3351569}). Moreover, the following holds:

\begin{Proposition}\label{ASuppAndLoc}
	For every object $M$ in $\cA$,
	\begin{equation*}
		\ASupp M=\setwithcondition{\alpha\in\ASpec\cA}{M_{\alpha}\neq 0}.
	\end{equation*}
\end{Proposition}

\begin{proof}
	This is proved in \cite[Proposition~6.2]{MR3351569} when $\cA$ is a Grothendieck category and the same proof works in the  general case.
\end{proof}

We state some basic properties on monoform objects. Recall that every nonzero subobject of a monoform object is again monoform (\cite[Proposition~2.2]{MR2964615}).

\begin{Proposition}\label{CharactOfMonoformObj}
	Let $H$ be an object in $\cA$ and $\alpha\in\ASpec\cA$. Then $H$ is a monoform object with $\overline{H}=\alpha$ if and only if
	\begin{enumerate}
		\item $H_{\alpha}$ is a simple object in $\cA_{\alpha}$, and
		\item $L_{\alpha}\neq 0$ for every nonzero subobject $L$ of $H$.
	\end{enumerate}
\end{Proposition}

\begin{proof}
	Assume that $H$ is a monoform object with $\overline{H}=\alpha$. For every nonzero subobject $L$ of $H$, we have $\alpha=\overline{H}=\overline{L}\in\ASupp L$, and hence $L_{\alpha}\neq 0$ by \cref{ASuppAndLoc}. Let $N'$ be a nonzero subobject of $H_{\alpha}$. Then by \cite[Corollary~4.3.10]{MR0340375}, there exists a subobject $N$ of $H$ satisfying $N_{\alpha}=N'$. By \cite[Proposition~2.14]{MR3351569}, we have $\alpha\notin\ASupp(H/N)$, and hence $H_{\alpha}/N_{\alpha}\cong (H/N)_{\alpha}=0$. This shows that $N'=N_{\alpha}=H_{\alpha}$, and $H_{\alpha}$ is simple.
	
	Assume that $H$ satisfies the latter two conditions. Let $L$ be a nonzero subobject of $H$. By the assumption, we have $(H/L)_{\alpha}\cong H_{\alpha}/L_{\alpha}=0$, and hence $N_{\alpha}=0$ for every subobject $N$ of $H/L$. Therefore no nonzero subobject of $H$ is isomorphic to a subobject of $H/L$. This shows that $H$ is monoform. Since we have $\alpha\in\ASupp H$, there exist subobjects $L_{1}\subset L_{2}\subset H$ such that $L_{2}/L_{1}$ is monoform and $\overline{L_{2}/L_{1}}=\alpha$. If $L_{1}\neq 0$, then $(L_{2}/L_{1})_{\alpha}\subset (H/L_{1})_{\alpha}=0$ by the above argument. This contradicts $\alpha=\overline{L_{2}/L_{1}}\in\ASupp(L_{2}/L_{1})$. Hence $L_{1}=0$, and $\overline{H}=\overline{L_{2}}=\alpha$.
\end{proof}

\begin{Lemma}\label{MorToMonoformObj}
	Let $\alpha\in\ASpec\cA$.
	\begin{enumerate}
		\item\label{MorFromTorsionObjToMonoformObj} Let $H$ be a monoform object in $\cA$ with $\overline{H}=\alpha$, and let $M$ be an object in $\cA$ satisfying $M_{\alpha}=0$. Then $\Hom_{\cA}(M,H)=0$.
		\item\label{MorBetweenMonoformObjs} Let $H_{1}$ and $H_{2}$ be monoform objects with $\overline{H_{1}}=\overline{H_{2}}=\alpha$. Then every nonzero morphism $H_{1}\to H_{2}$ is a monomorphism.
	\end{enumerate}
\end{Lemma}

\begin{proof}
	\cref{MorFromTorsionObjToMonoformObj} Assume that a morphism $f\colon M\to H$ is nonzero. Then we have $\alpha=\overline{H}=\overline{\Im f}\in\ASupp(\Im f)$. On the other hand, since $\Im f$ is a quotient object of $M$, we have $(\Im f)_{\alpha}=0$ by the assumption. By \cref{ASuppAndLoc}, we have $\alpha\notin\ASupp(\Im f)$. This is a contradiction.
	
	\cref{MorBetweenMonoformObjs} Assume that a morphism $f\colon H_{1}\to H_{2}$ is not a monomorphism. Then \cref{CharactOfMonoformObj} shows that $(\Ker f)_{\alpha}$ is a nonzero subobject of the simple object $(H_{1})_{\alpha}$, and we have $(H_{1}/\Ker f)_{\alpha}\cong (H_{1})_{\alpha}/(\Ker f)_{\alpha}=0$. By \cref{MorFromTorsionObjToMonoformObj}, the induced monomorphism $H_{1}/\Ker f\to H_{2}$ is zero, and we obtain $f=0$.
\end{proof}

The following lemma will be used in the next section.

\begin{Lemma}\label{SubobjOfQuotCat}
	Let $\cA$ be a Grothendieck category. Let $M$ be an object in $\cA$ and $\alpha\in\ASpec\cA$. For each subobject $L'$ of $M_{\alpha}$, there exists a subobject $L$ of $M$ such that
	\begin{enumerate}
		\item $L_{\alpha}=L'$,
		\item $L$ is largest with respect to the property $L_{\alpha}\subset L'$, and
		\item $N_{\alpha}\neq 0$ for every nonzero subobject $N$ of $M/L$.
	\end{enumerate}
\end{Lemma}

\begin{proof}
	Let $L$ be the sum of all subobjects $N$ of $M$ satisfying $N_{\alpha}\subset L'$. Since the equality is satisfied for some $N$ by \cite[Corollary~4.3.10]{MR0340375}, $L_{\alpha}=L'$. The other conditions are also satisfied.
\end{proof}

The notion of a compressible object plays an important role in the proof of the main result. We recall the definition and some basic properties of compressible objects.

\begin{Definition}\label{CompressibleObj}
	A nonzero object $H$ in $\cA$ is called \emph{compressible} if each nonzero subobject $L$ of $H$ has some subobject that is isomorphic to $H$.
\end{Definition}

\begin{Proposition}\label{PropertiesOfCompressibleObj}\leavevmode
	\begin{enumerate}
		\item\label{SubOfCompressibleIsCompressible} Every nonzero subobject of a compressible object in $\cA$ is compressible.
		\item\label{CompressibleObjInLocNoethGrothCatIsNoethAndMonoform} If $\cA$ has a generating set consisting of noetherian objects, then every compressible object in $\cA$ is noetherian and monoform.
	\end{enumerate}
\end{Proposition}

\begin{proof}
	\cref{SubOfCompressibleIsCompressible} Straightforward.
	
	\cref{CompressibleObjInLocNoethGrothCatIsNoethAndMonoform} Let $H$ be a compressible object in $\cA$. Then $H$ has a nonzero noetherian subobject, and by \cite[Theorem~2.9]{MR2964615}, the subobject has a monoform subobject $L$. There exists a subobject of $L$ that is isomorphic to $H$, and hence $H$ is monoform.
\end{proof}

In the case of commutative rings, compressible modules are characterized as follows.

\begin{Proposition}\label{CompressibleObjForCommRing}
	Let $R$ be a commutative ring.
	\begin{enumerate}
		\item\label{CompressibleObjAndPrimeIdeal} \textnormal{(\cite[3.38]{MR0393118})} For an ideal $\kp$ of $R$, the object $R/\kp$ in $\Mod R$ is compressible if and only if $\kp$ is a prime ideal.
		\item\label{CharactOfCompressibleObjForCommRing} An $R$-module $H$ is a compressible object in $\Mod R$ if and only if $H$ is isomorphic to a nonzero $R$-submodule of $R/\kp$ for some $\kp\in\Spec R$.
	\end{enumerate}
\end{Proposition}

\begin{proof}
	\cref{CompressibleObjAndPrimeIdeal} For the convenience of the reader, we include a proof.
	
	Assume that $R/\kp$ is compressible. Then for each $a\in R\setminus\kp$, we have $\Ann_{R}(\overline{a}\in R/\kp)=\kp$ since $\overline{a}R$ has an $R$-submodule isomorphic to $R/\kp$. This shows that $\kp$ is a prime ideal.
	
	Conversely, let $\kp$ is a prime ideal, and let $L$ be a nonzero $R$-submodule of $R/\kp$. Take a nonzero element $\overline{a}\in L$. Since we have $\Ann_{R}(\overline{a})=\kp$, it holds that $\overline{a}R\cong R/\kp$.
	
	\cref{CharactOfCompressibleObjForCommRing} Assume that $H$ is compressible. Take a nonzero element $x\in H$. Then $H$ is isomorphic to an $R$-submodule of $xR$, which is compressible by \cref{PropertiesOfCompressibleObj} \cref{SubOfCompressibleIsCompressible}. We have $xR\cong R/\Ann_{R}(x)$, and by \cref{CompressibleObjAndPrimeIdeal}, the ideal $\Ann_{R}(x)$ is prime.
	
	The converse follows from \cref{CompressibleObjAndPrimeIdeal} and \cref{PropertiesOfCompressibleObj} \cref{SubOfCompressibleIsCompressible}.
\end{proof}

\begin{Remark}\label{CompressibleObjForNoncommRing}
	For a commutative ring $R$, each atom in $\Mod R$ is of the form $\overline{R/\kp}$ for some $\kp\in\Spec R$, and $\overline{R/\kp}$ is a compressible object in $\Mod R$.
	
	On the other hand, for a right noetherian ring, an atom is not necessarily represented by a compressible object. Indeed, Goodearl \cite{MR588447} constructed a left and right noetherian ring $\Lambda$ such that there exists a monoform $\Lambda$-module that does not have a compressible $\Lambda$-submodule.
\end{Remark}

%%%%%%%%%%%%%%%%%%%%%%%%%%%%%%%%%%%%%%%%%%%%%%%%%%%%%%%%%%%%%%%%%%%%%%%%%%%%%%%%
\section{Proof of the main result}
\label{sec.ProofOfMainResults}
%%%%%%%%%%%%%%%%%%%%%%%%%%%%%%%%%%%%%%%%%%%%%%%%%%%%%%%%%%%%%%%%%%%%%%%%%%%%%%%%

Throughout this section, let $\cA$ be a Grothendieck category having a noetherian generator unless otherwise specified. The following theorem shows that the behavior of minimal atoms is similar to atoms in $\Mod R$ for a commutative ring $R$.

\begin{Theorem}\label{MinimalAtomIsRepresentedByCompressibleObj}
	For every $\alpha\in\AMin\cA$, there exists a compressible object $H$ in $\cA$ satisfying $\alpha=\overline{H}$.
\end{Theorem}

\begin{proof}
	Let $U$ be a noetherian generator in $\cA$. By \cite[Proposition~6.6 (1)]{MR3351569}, the topological space $\ASpec\cA_{\alpha}$ is homeomorphic to the singleton $\set{\alpha}$ since $\alpha$ is minimal. Hence by \cref{ImageOfGeneratorAndNoethObj} \cref{ImageOfNoethObj} and \cite[Proposition~3.7 (3)]{MR3351569}, the object $U_{\alpha}$ in $\cA_{\alpha}$ has a composition series
	\begin{equation*}
		0=L'_{0}\subset L'_{1}\subset\cdots\subset L'_{n}=U_{\alpha}.
	\end{equation*}
	By using \cref{SubobjOfQuotCat}, we obtain a filtration
	\begin{equation*}
		L_{0}\subset L_{1}\subset\cdots\subset L_{n}=U
	\end{equation*}
	such that $(L_{i})_{\alpha}=L'_{i}$ for each $i=0,\ldots,n$, and $N_{\alpha}\neq 0$ for every nonzero subobject $N$ of $U/L_{i}$. For each $i=1,\ldots,n$, the object $(L_{i}/L_{i-1})_{\alpha}\cong L'_{i}/L'_{i-1}$ is simple, and $N_{\alpha}\neq 0$ for every nonzero subobject $N$ of $L_{i}/L_{i-1}$. Hence by \cref{CharactOfMonoformObj}, the object $H_{i}:=L_{i}/L_{i-1}$ is monoform and $\overline{H_{i}}=\alpha$. Since the monoform objects $H_{1},\ldots,H_{n}$ represent the same atom $\alpha$, there exists a monoform object $H$ in $\cA$ such that for each $i=1,\ldots,n$, the object $H$ is isomorphic to some subobject of $H_{i}$.
	
	Let $H'$ be a nonzero subobject of $H$. Since $U$ is a generator in $\cA$, there exists a nonzero morphism $U\to H'$. By using the filtration of $U$, we obtain a nonzero morphism $L_{0}\to H'$ or a nonzero morphism $L_{i}/L_{i-1}\to H'$ for some $i=1,\ldots,n$. \cref{MorToMonoformObj} \cref{MorFromTorsionObjToMonoformObj} implies that the former case does not occur, and \cref{MorToMonoformObj} \cref{MorBetweenMonoformObjs} implies that the nonzero morphism $H_{i}=L_{i}/L_{i-1}\to H'$ is a monomorphism. Thus $H'$ has a subobject that is isomorphic to $H$ by the definition of $H$. This shows that $H$ is compressible.
\end{proof}

The following result shows the existence of minimal atoms, which will be used in the proof of \cref{AMinIsDiscreteAndClosed}.

\begin{Proposition}\label{ExistenceOfMinimalAtom}
	For each $\alpha\in\ASpec\cA$, there exists $\beta\in\AMin\cA$ satisfying $\beta\leq\alpha$.
\end{Proposition}

\begin{proof}
	Let $U$ be a noetherian generator in $\cA$. By \cref{ASuppOfGenerator}, $\ASpec\cA=\ASupp U$. Thus the statement follows from \cite[Proposition~4.7]{MR3351569}.
\end{proof}

In order to prove the finiteness of $\AMin\cA$, we show some topological properties of it. Although the statement is quite simple, the proof needs our previous results.

\begin{Proposition}\label{AMinIsDiscreteAndClosed}
	$\AMin\cA$ is a closed subset of $\ASpec\cA$, and the induced topology on $\AMin\cA$ is discrete.
\end{Proposition}

\begin{proof}
	Let $\alpha\in\AMin\cA$. We show that
	\begin{equation*}
		V(\alpha):=\setwithcondition{\beta\in\ASpec\cA}{\alpha\leq\beta}\ \ \text{and}\ \ W(\alpha):=\setwithcondition{\beta\in\ASpec\cA}{\alpha<\beta}
	\end{equation*}
	are open subsets of $\ASpec\cA$.
	
	By \cref{MinimalAtomIsRepresentedByCompressibleObj}, there exists a compressible object $H$ in $\cA$ satisfying $\alpha=\overline{H}$. It follows from \cite[Proposition~4.2]{MR3351569} that $V(\alpha)=\ASupp H$, which is an open subset of $\ASpec\cA$. Since the minimal atom $\alpha$ is a closed point in $\ASpec\cA$, the subset $W(\alpha)=V(\alpha)\setminus\set{\alpha}$ of $\ASpec\cA$ is open.
	
	By \cref{ExistenceOfMinimalAtom}, we obtain the equation
	\begin{equation*}
		\AMin\cA=\ASpec\cA\setminus\bigcup_{\alpha\in\AMin\cA}W(\alpha),
	\end{equation*}
	which implies that $\AMin\cA$ is a closed subset of $\ASpec\cA$. The equation
	\begin{equation*}
		V(\alpha)\cap\AMin\cA=\set{\alpha}
	\end{equation*}
	shows that $\alpha$ is an open point in $\AMin\cA$.
\end{proof}

We are ready to prove the main result of this paper.

\begin{Theorem}\label{AMinOfGrothCatWithNoethGenIsFinite}
	Let $\cA$ be a Grothendieck category having a noetherian generator. Then $\AMin\cA$ is a finite set.
\end{Theorem}

\begin{proof}
	By \cref{AMinIsDiscreteAndClosed}, the subset $\Phi:=\ASpec\cA\setminus\AMin\cA$ of $\ASpec\cA$ is open. Let $\cX:=\ASupp^{-1}\Phi$. Then by \cite[Theorem~5.17]{MR3351569} and \cref{ClassifOfLocSubcatAndSerreSubcat} \cref{ClassifOfLocSubcat}, the topological space $\ASpec(\cA/\cX)$ is homeomorphic to $\AMin\cA$. Since it is discrete, by \cref{ImageOfGeneratorAndNoethObj} and \cite[Proposition~3.7 (3)]{MR3351569}, the image $U'$ of a noetherian generator in $\cA$ by the canonical functor $\cA\to\cA/\cX$ is a generator of finite length in $\cA/\cX$. Let $S'_{1},\ldots,S'_{n}$ be the composition factors of $U'$. It follows from \cref{ASuppOfGenerator} that
	\begin{equation*}
		\ASpec(\cA/\cX)=\ASupp U'=\ASupp S'_{1}\cup\cdots\cup\ASupp S'_{n}=\set{\overline{S'_{1}},\ldots,\overline{S'_{n}}}.\qedhere
	\end{equation*}
	Therefore $\AMin\cA$ is also a finite set.
\end{proof}

It is known that there is a bijective correspondence between the equivalence classes of locally noetherian Grothendieck categories and those of noetherian abelian categories (\cite[Theorem~II.4.1]{MR0232821}; see also \cite[section~5.8]{MR0340375}). For a locally noetherian Grothendieck category $\cA$, the corresponding noetherian abelian category is the Serre subcategory $\noeth\cA$ of $\cA$ consisting of all noetherian objects. It is easy to see that $\cA$ has a noetherian generator if and only if $\noeth\cA$ has a generator. This leads us to the following corollary:

\begin{Corollary}\label{AMinOfNoethAbelCatWithGenIsFinite}
	Let $\cA$ be a noetherian abelian category having a generator. Then $\AMin\cA$ is a finite set.
\end{Corollary}

\begin{proof}
	By the correspondence described above, there exists a Grothendieck category $\cG$ having a noetherian generator such that $\noeth\cG\cong\cA$. By \cite[Proposition~5.3]{MR2964615}, $\ASpec\cG$ is homeomorphic to $\ASpec(\noeth\cG)$, and hence to $\ASpec\cA$. Therefore $\AMin\cA$ bijectively corresponds to $\AMin\cG$, which is a finite set by \cref{AMinOfGrothCatWithNoethGenIsFinite}.
\end{proof}

The main results can be applied to a right noetherian ring. Note that $\noeth(\Mod\Lambda)$ is the full subcategory $\mod\Lambda$ consisting of all finitely generated right $\Lambda$-modules.

\begin{Corollary}\label{FinitenessOfAMinForOneSidedNoethRing}
	Let $\Lambda$ be a right noetherian ring. Then the set $\AMin(\Mod\Lambda)\cong\AMin(\mod\Lambda)$ is a finite set.
\end{Corollary}

\begin{proof}
	This follows from \cref{AMinOfGrothCatWithNoethGenIsFinite} since $\Lambda$ is a noetherian generator in $\Mod\Lambda$.
\end{proof}

\begin{Remark}\label{CompactnessOfASpec}
	Let $\cA$ be a Grothendieck category having a noetherian generator (or a noetherian abelian category having a generator). Then \cite[Proposition~3.6]{MR3351569} and the argument in the proof of \cref{ExistenceOfMinimalAtom} show that $\ASpec\cA$ is compact with respect to the localizing topology. Since every open subset of $\ASpec\cA$ that contains a minimal atom $\alpha$ also contains all atoms $\beta$ with $\alpha\leq\beta$, \cref{AMinOfGrothCatWithNoethGenIsFinite} also explains why $\ASpec\cA$ is compact.
\end{Remark}

\begin{Remark}\label{OtherProofForModCat}
	\cref{FinitenessOfAMinForOneSidedNoethRing} can also be deduced from a result of Beachy \cite[Theorem~3.6]{MR0327813}: There is a bijection between the \emph{maximal torsion radicals} of $\Mod\Lambda$ and the minimal prime two-sided ideals of $\Lambda$ (see \cite[section~1]{MR0327813} for the definition of (maximal) torsion radicals).
	
	For a Grothendieck category $\cA$, there exists an order-preserving bijection between the torsion radicals of $\cA$ and the localizing subcategories of $\cA$ (see \cite[Proposition~VI.3.1]{MR0389953}; note that the terminologies is different there). If $\cA$ is locally noetherian, then the localizing subcategories bijectively correspond to the open subsets of $\ASpec\cA$ by \cref{ClassifOfLocSubcatAndSerreSubcat} \cref{ClassifOfLocSubcat}, and also to the closed subsets of $\ASpec\cA$ by taking complements. Hence the maximal torsion radicals of $\cA$ bijectively correspond to the minimal nonempty closed subsets of $\cA$. By \cref{ExistenceOfMinimalAtom} together with the definition of the specialization order on $\ASpec\cA$, a subset of $\ASpec\cA$ is a minimal nonempty closed subset if and only if it consists of exactly one minimal atom. Therefore the maximal torsion radicals of $\cA$ bijectively corresponds to $\AMin\cA$. This applies to the case $\cA=\Mod\Lambda$.
	
	On the other hand, it is well known that a right noetherian ring has only finitely many minimal prime two-sided ideals (see \cite[Theorem~3.4]{MR2080008}). This completes another proof of \cref{FinitenessOfAMinForOneSidedNoethRing} for $\Mod\Lambda$.
\end{Remark}

\begin{Example}\label{ASpecOfDomain}
	Let $\Lambda$ be a right noetherian ring that is a domain, that is, for each $a,b\in\Lambda$, the condition $ab=0$ implies $a=0$ or $b=0$. Then the object $\Lambda$ in $\Mod\Lambda$ is compressible since the right $\Lambda$-module $a\Lambda$ is isomorphic to $\Lambda$ for each $a\in\Lambda$. Hence by \cref{PropertiesOfCompressibleObj} \cref{CompressibleObjInLocNoethGrothCatIsNoethAndMonoform}, the object $\Lambda$ in $\Mod\Lambda$ is monoform. Since $\Lambda$ is a generator in $\Mod\Lambda$, it follows from \cref{ASuppOfGenerator} that $\ASpec(\Mod\Lambda)=\ASupp\Lambda$. We deduce from \cite[Proposition~4.2]{MR3351569} that $\AMin(\Mod\Lambda)=\set{\overline{\Lambda}}$. From the viewpoint of \cref{OtherProofForModCat}, the unique minimal atom $\overline{\Lambda}$ corresponds to the zero ideal of $\Lambda$.
\end{Example}

Wu \cite{MR1029695} gave an example of a Grothendieck category having a noetherian generator that does not satisfy the condition Ab4*, that is, direct products are not exact. In particular, the Grothendieck category is not equivalent to the category of right modules over a ring. We give a description of the minimal atoms of such Grothendieck categories.

\begin{Example}\label{ExOfGrothCatWithNoethGenerator}
	Let $R$ be a commutative noetherian Gorenstein ring with Krull dimension at least $2$. Let $\Phi:=\setwithcondition{\kp\in\Spec R}{\height\kp\geq 2}$ and
	\begin{equation*}
		\cX:=\Supp^{-1}\Phi=\setwithcondition{M\in\Mod R}{\Supp M\subset\Phi}.
	\end{equation*}
	Wu \cite{MR1029695} showed (in a more general setting) that the Grothendieck category $\cA:=(\Mod R)/\cX$ does not satisfy Ab4*, whereas the image of $R$ by the canonical functor $\Mod R\to\cA$ is a noetherian generator in $\cA$.
	
	In this case, the atom spectrum $\ASpec\cA$ is isomorphic to $(\Spec R)\setminus\Phi=\setwithcondition{\kp\in\Spec R}{\height\kp\leq 1}$ as a partially ordered set. Therefore $\AMin\cA$ is in bijection with the set of minimal prime ideals of $R$.
\end{Example}

\begin{Remark}\label{ProofForCommRing}
	For a commutative noetherian ring $R$, the finiteness of the number of minimal prime ideals can be shown by using the equation
	\begin{equation*}
		\Min(\Spec R)=\Min(\Supp R)=\Min(\Ass R)
	\end{equation*}
	and the fact that $\Ass R$ is a finite set. For a Grothendieck category $\cA$ having a noetherian generator $U$, it still holds that $\Min(\ASpec\cA)=\Min(\ASupp U)$ by \cref{ASuppOfGenerator}, and the set $\AAss U$ of \emph{associated atoms} of $U$ is finite (see \cite[section~3]{MR2964615} for the definition and basic properties of associated atoms). However, the equality $\Min(\ASupp U)=\Min(\AAss U)$ is far from being true, even in the case where $\cA=\Mod\Lambda$ for a ring $\Lambda$. For example, let $\Lambda$ be the ring
	\begin{equation*}
		\begin{bmatrix}
			K & 0 \\
			K & K
		\end{bmatrix}
	\end{equation*}
	of $2\times 2$ lower triangular matrices over a field $K$. Define the indecomposable right $\Lambda$-modules
	\begin{equation*}
		S_{1}:=
		\begin{bmatrix}
			K & 0
		\end{bmatrix}
		,\ S_{2}:=
		\begin{bmatrix}
			K & K
		\end{bmatrix}
		/
		\begin{bmatrix}
			K & 0
		\end{bmatrix}
		,\ P_{2}:=
		\begin{bmatrix}
			K & K
		\end{bmatrix}.
	\end{equation*}
	Then $\Lambda=S_{1}\oplus P_{2}$ as a right $\Lambda$-module, and we have
	\begin{equation*}
		\AMin(\Mod\Lambda)=\Min(\ASupp\Lambda)=\ASupp\Lambda=\set{\overline{S_{1}},\,\overline{S_{2}}}\neq\set{\overline{S_{1}}}=\AAss\Lambda=\Min(\AAss\Lambda).
	\end{equation*}
\end{Remark}

%%%%%%%%%%%%%%%%%%%%%%%%%%%%%%%%%%%%%%%%%%%%%%%%%%%%%%%%%%%%%%%%%%%%%%%%%%%%%%%%

%%%%%%%%%%%%%%%%%%%%%%%%%%%%%%%%%%%%%%%%%%%%%%%%%%%%%%%%%%%%%%%%%%%%%%%%%%%%%%%%

\end{document}